\documentclass[12pt]{amsart}

\usepackage{amsmath}
\usepackage{amssymb}
\usepackage{amsthm}
\usepackage{amscd}
\usepackage[all]{xy}
\usepackage[utf8]{inputenc}
\usepackage[T1]{fontenc}   
\usepackage{enumerate}
\usepackage{color}
\usepackage{url}
\usepackage{hyperref}

\usepackage{mathpazo}
\usepackage[scaled=.95]{helvet}
\usepackage{courier}

\swapnumbers
\headsep=1cm \numberwithin{equation}{section}
\newtheorem{theorem}{Theorem}[section]
\newtheorem{lemma}[theorem]{Lemma}
\newtheorem{proposition}[theorem]{Proposition}
\newtheorem{corollary}[theorem]{Corollary}

\theoremstyle{definition}
\newtheorem{definition}[theorem]{Definition}

\newtheorem{remark}[theorem]{Remark}
\newtheorem{remark and definition}[theorem]{Remark and Definition}
\newtheorem{remark and notation}[theorem]{Remark and Notation}

\newtheorem{conjecture}[theorem]{Conjecture}
\newtheorem{question}[theorem]{Question}

\newtheorem*{acknowledgement}{Acknowledgement}

\newtheorem{facts}[theorem]{Facts}
\newtheorem{questions}[theorem]{Questions}

\newcommand\Spec{\operatorname{Spec}}
\newcommand\Hom{\operatorname{Hom}}
\newcommand\Ext{\operatorname{Ext}}

\newcommand\depth{\operatorname{depth}}

\newcommand\Ker{\operatorname{\Ker}}

\newcommand\pd{\operatorname{pd}}

\newcommand\Supp{\operatorname{Supp}}

\author[Mendoza-Rubio,\,  Jorge-P\'erez]{Victor D. Mendoza-Rubio and Victor H. Jorge-P\'erez}

\title[On modules whose dual is of finite Gorenstein dimension]{On modules whose dual is of finite Gorenstein dimension}

\address{Universidade de S{\~a}o Paulo -
ICMC, Caixa Postal 668, 13560-970, S{\~a}o Carlos-SP, Brazil}
\email{vicdamenru@usp.br}

\address{Universidade de S{\~a}o Paulo -
ICMC, Caixa Postal 668, 13560-970, S{\~a}o Carlos-SP, Brazil}
\email{vhjperez@icmc.usp.br}

\keywords{Gorenstein dimension, Algebraic dual, $k$-torsion, Totally reflexive modules, Free modules, Vanishing of Ext, Auslander-Reiten conjecture, Gorenstein ring, Module of K\"ahler differentials, Herzog-Vasconcelos conjecture}
\subjclass[2020]{Primary: 13D45, 13D07, 13C10, 13C14; Secondary: 13D05, 13D02, 13H10, 14B15, 14B05.}

\begin{document}

\begin{abstract} 
In this paper, we aim to obtain some results under the condition that the dual of a module over a commutative Noetherian ring has finite Gorenstein dimension. In this direction, we derive results involving vanishing of Ext as well as the freeness or totally reflexivity of modules. For instance,  we provide a generalization of a celebrated theorem by Auslander and Bridger, obtain criteria for the  totally reflexivity of modules over Cohen-Macaulay rings as well as of locally totally reflexive  modules on the punctured spectrum, and recover a result by Araya. Moreover, we prove that the Auslander-Reiten conjecture holds true for all finitely generated modules $M$ over a commutative Noetherian ring $R$ such that $\operatorname{G-dim}_R(\Hom_R(M,R))<\infty$ and $\operatorname{pd}_R(\Hom_R(M,M))<\infty$. Additionally, we derive Gorenstein criteria under the condition that the dual of certain modules is of finite Gorenstein dimension. Furthermore,
we explore some applications in the theory of the modules of K\"ahler differentials of order $n\geq 1$, specifically concerning the $k$-torsionfreeness of these modules and the Herzog-Vasconcelos conjecture. 
\end{abstract}

\maketitle
\section{Introduction}
The Gorenstein dimension is a well-developed topic in commutative algebra, as evidenced by various works such as \cite{Stablemoduletheory, LinkageOfModulesAndTheSerreConditions, ktorsionlessmoduleswithfiniteGorensteindimension, GorensteinDimensionAndTorsionOfModulesOverCommutativeNoetherianRings,  GorensteinringsviahomologicaldimensionsandsymmetryinvanishingofExtandTatecohomology, HomologicalDimensionsOfRigidModules}. Despite the extensive research in this area, there are relatively few results concerning the condition of the Gorenstein dimension of the dual of certain modules being finite. The aim of this work is to obtain results in this direction, including some that involve vanishing of Ext and homological criteria.

 Let $R$ be a (commutative) Noetherian ring.  We prove the following theorem, which is a generalization of a result that  Evans and Griffith attribute to Auslander and Bridger in \cite[Theorem 3.8]{Syzygies} (see also \cite[Lemma 1.3]{AscentOfFiniteCohen–MacaulayType}) and is an analogous version of \cite[Proposition 2.4]{LinkageOfModulesAndTheSerreConditions} when $C=R$. \begin{theorem}[{See Theorem \ref{GteosPa}}] \label{IGteosPa}
    Let $k$ be a non-negative integer. Let $M$ be a finitely generated $R$-module such that $M^\ast$ has locally finite Gorenstein dimension on $\widetilde{X}^{k-1}(R)$. Then: 
    \begin{center}
        $M$ is $k$-torsionfree $\Longleftrightarrow$ $M$ is $k$-syzygy $\Longleftrightarrow$ $M$ satisfies $(\widetilde{S}_k)$.
    \end{center}
\end{theorem}
This theorem provides an answer to \cite[Question 1.1]{WhenArenSysymodulesnTorsionfree}. Additionally, it allows us to recover one of the results described in the main theorem of \cite{HomologicalDimensionsOfRigidModules}.
\begin{corollary}[{\cite[Theorem 5.8(1)]{HomologicalDimensionsOfRigidModules}}] \label{I516s} Suppose that $R$ is local. 
    Let $M$ be a finitely generated $R$-module such that $\operatorname{G-dim}_R(M^\ast)<\infty$ and $k$ be a non-negative integer. If $M$ satisfies $(\widetilde{S}_k)$, then $M$ is $k$-torsionfree.
\end{corollary}

As an application of Corollary \ref{I516s} and \cite[
Theorem 42]{GorensteinDimensionAndTorsionOfModulesOverCommutativeNoetherianRings}, we characterize the totally reflexivity of a module whose dual is of finite Gorenstein dimension in terms of certain Ext vanishing and some other conditions. We obtain the following. 
\begin{theorem}[See Theorems \ref{pos5}, \ref{genCM} and \ref{seaa1}] \label{Is4ap6} 
    Suppose that $R$ is local of depth $t$. Let $M$ be a finitely generated $R$-module such that $\operatorname{G-dim}_R(M^\ast)<\infty$. Then:
    \begin{enumerate}
        \item $M$ satisfies $(\widetilde{S}_t)$ if and only if $M$ is totally reflexive if and only if $\Ext_R^i(M,R)=0$ for all $1\leq i \leq t$.
        \item Let $R$ be a Cohen-Macaulay local ring. If $n$ is a non-negative integer such that $\operatorname{depth}_R(M)\geq n$ and $\Ext_R^i(M,R)=0$ for all $1\leq i \leq t-n,$ then $M$ is totally reflexive. 
        \item If $M$ is non-zero and locally totally reflexive on the punctured spectrum of $R$, then $\operatorname{depth}_R(M)\leq t$ with equality if and only if $M$ is totally reflexive. 
    \end{enumerate}
\end{theorem}

Additionally, part of this work is motivated by the celebrated Auslander-Reiten conjecture (\cite{OnageneralizedversionoftheNakayamaconjecture}), which states the following. 
\begin{conjecture}[Auslander-Reiten] \label{ARC} Let $M$ be a finitely generated $R$-module. If $\Ext_R^i(M,R)=\Ext_R^i(M,M)=0$ for all $i>0$, then $M$ is projective. 
\end{conjecture}

This conjecture is still open but it is true in several cases.  For a list of some of them, we refer the reader to the introductions of \cite{injectivedimensiontakahashi, AuslanderReitenConjectureforNormalRings,FiniteHomologicalDimensionOfHomv3}. In \cite{TheAuslander-ReitenconjectureforGorensteinrings}, Araya proved the following result that provided a positive answer to Auslander-Reiten conjecture, and that implies that it is valid for normal Gorenstein local rings of dimension at least two.

\begin{theorem}[{\cite[Corollary 10]{TheAuslander-ReitenconjectureforGorensteinrings}}] \label{araq}
Suppose that $R$ is Gorenstein local of dimension $d\geq 2$. Let $M$ be a maximal Cohen-Macaulay $R$-module. If $M$ is locally free on the punctured spectrum of $R$ and $\Ext_R^{d-1}(M,M)=0$, then $M$ is free.
\end{theorem}
Subsequently, results concerning the context of Gorenstein rings that recover this theorem have appeared. For instance, the last theorem in \cite[Section 1]{AnAuslander–Reitenprincipleinderivedcategories} and \cite[Corollary 1.6]{OnthevanishingofselfextensionsoverCohen-Macaulaylocalrings}. It is important to mention that in a more general context, \cite[Proposition 2.10(4)]{AuslanderReitenConjectureforNormalRings} and \cite[Proposition 1.7]{TwoTheoremsOnTheVanishingOfExt} are improvements of these results. In terms of the finiteness of the Gorenstein dimension of the dual of a module, we achieve the following results that recover the last theorem in \cite[Section 1]{AnAuslander–Reitenprincipleinderivedcategories},  \cite[Corollary 1.6]{OnthevanishingofselfextensionsoverCohen-Macaulaylocalrings} and, of course,  Theorem \ref{araq}. 
\begin{theorem}[See Corollary \ref{432c}]  Suppose that $R$ is local. Let  $M$ be a finitely generated $R$-module such that $\operatorname{G-dim}_R(M^\ast)<\infty$. Suppose that $M$ is locally free on the punctured spectrum of $R$. If $\depth(M)\geq t$ and $\Ext_R^{t-1}(M,M)=0$, then $M$ is free. 
\end{theorem}

\begin{theorem}[See Theorem \ref{rpsCM}]\label{IrpsCM}
    Let $R$ satisfy $(S_1)$ and  let $X$ be a subset of $\Spec(R)$ containing $X^1(R)$. Let 
    \begin{center}
        $s:=\inf \{\operatorname{depth}_{R_\mathfrak{p}} (R_{\mathfrak{p}}) \mid \mathfrak{p} \in \operatorname{Spec}(R) \backslash X\}$ and $t:=\sup \{\operatorname{depth}_{R_\mathfrak{p}} (R_{\mathfrak{p}}) \mid \mathfrak{p} \in \operatorname{Spec}(R) \backslash X\}$.
    \end{center}
    Let $M$ be a finitely generated $R$-module. Suppose that $M$ is locally free on $X$ and that $\operatorname{G-dim}_R(M^\ast)<\infty$. If $\Ext_R^i(M,R)=\Ext_R^j(M,M)=0$ for all $1\leq i \leq t$ and $s-1 \leq j \leq t-1$, then $M$ is projective.
\end{theorem}

On the other hand, in \cite[Corollary 6.9(2)]{FiniteHomologicalDimensionOfHomv3}, Dey and Ghosh proved that the Auslander-Reiten conjecture holds for finitely generated modules $M$ over a (commutative) Noetherian ring $R$ such that $\operatorname{G-dim}_R(M)<\infty$ and $\operatorname{pd}_R(\Hom_R(M,M))<\infty$. Motivated by this, we analyze if this fact holds if we consider $\operatorname{G-dim}_R(M^\ast)<\infty$ instead of $\operatorname{G-dim}_R(M)<\infty$. Consequently, we get the following result, which not only affirms that this holds but also strengthens \cite[Corollary 2.14]{injectivedimensiontakahashi}. 
\begin{theorem}[See Theorem \ref{fjk}] \label{Ifjk}
     Suppose that $R$ is local of depth $t$. Let $M$  be a finitely generated $R$-module such that $\operatorname{G-dim}_R(M^\ast)<\infty$. Then the following conditions are equivalent:
     \begin{enumerate}[(1)]
         \item $M$ is free.
         \item $\Hom_R(M,M)$ is free and $\Ext_R^j(M,M)=0$ for all $1\leq j \leq t-1$.
         \item  $\Hom_R(M,M)$ has finite projective dimension, and $\Ext_R^i(M,R)= \Ext_R^j(M,M) =0$ for all $1\leq i \leq t$ and $1\leq j \leq t-1$. 
     \end{enumerate}
\end{theorem}

Additionally, we also derive some Gorenstein criteria involving the condition that the dual of certain modules are of finite Gorenstein dimension. 

\begin{theorem}[See Propositions \ref{Gorentein2}, \ref{sq884} and \ref{Gdual}, and  Theorem \ref{GG}] \label{GorC} Suppose that $(R, \mathfrak{m})$ is local. Then $R$ is Gorenstein in each one of the following cases:
\begin{enumerate}
    \item  $R$ is a Cohen-Macaulay local ring, and there exists an $R$-module $M$ such that \linebreak $\operatorname{G-dim}_R(M^\ast)<\infty$ and  that satisfies any one of the following:
    \begin{enumerate}[(1.1)]
        \item $M$ is Cohen-Macaulay of positive rank and $2\mu(M)>e(R)\operatorname{rank}(M)$.
        \item $M$ is an Ulrich module and $M^\ast \not=0$. 
    \end{enumerate}
    \item $\operatorname{G-dim}_R(M^\ast)<\infty$ for all finitely generated $R$-modules $M$. 
    \item  $R$ has positive depth, and every finitely generated $R$-module $M$ whose dual $M^\ast$ is non-zero and of finite Gorenstein dimension, is also of finite Gorenstein dimension. 
    \item  $R$ has depth at most one and  $\operatorname{G-dim}_R(\mathfrak{m}^\ast)<\infty$. 
    
    \end{enumerate}
\end{theorem}

In the last part of this paper, as an application, we focus on the Kähler differential modules $\Omega_{X/Y}^{(n)}$ and the derivation modules $\mathcal{D}er_Y^{n}(X)$, defined over locally Noetherian schemes, where $n \geq 1$ is an integer. Initially, we investigate the conditions under which the Kähler differential module $\Omega_{X/Y}^{(n)}$ is  reflexive, $k$-torsionfree, and $k$-syzygy for $k=1,2$ (see Theorem \ref{stalksquema1}, and Corollaries \ref{corsche1} and \ref{corsche2}). The analysis of these properties in Kähler differential modules of order $n$, in the context of regularity questions of affine rings, has a long history in the literature, being discussed by several authors, especially in the context of algebraic varieties and analytic spaces, particularly when $Y$ is the spectrum of a field (see, for example, \cite{Kunz, Milher, graf, Lipman, Suzuki, Vetter, Lebelt, Greuel}).

Furthermore, we consider the Herzog-Vasconcelos conjecture, which states: Let $R$ be a local $k$-algebra, where $k$ is a field of characteristic zero. If ${\rm pd}_R({\rm Der}_k(R)) < \infty$, then ${\rm Der}_k(R)$ is a free $R$-module. This conjecture remains open, although some particular cases have already been resolved (see \cite[Section 4]{He-M}, \cite{HV, PerezNeto}). Motivated by this conjecture, we propose a generalized version in the context of schemes, called the \textit{generalization of the Herzog-Vasconcelos conjecture} (GHVC), which is described as follows: Let $n \geq 1$ and ${\rm pd}_{\mathcal{O}_X}{\mathcal{D}er}_Y^n(X) < \infty$. Under what conditions on $X$ and $Y$, and for which values of $n$, does this imply that ${\mathcal{D}er}_Y^n(X)$ is locally free? (see Question \ref{questionHV}). Following this line of investigation, we present an application, in which we provide a partial answer to Question \ref{questionHV} (see Proposition \ref{diff1}).

The organization of this paper is as follows. In Section \ref{section2} we provide definitions, notations, and some results that are considered in this paper. In Section 3 we prove Theorem \ref{IGteosPa} and explore some of its consequences. In Section 4, we prove  Theorem \ref{Is4ap6}, and in Section 5 we provide freeness criteria including  Theorems \ref{IrpsCM} and \ref{Ifjk}. In Section 6, we discuss some questions about the Gorensteiness of a ring involving Gorenstein dimension of the dual of modules, and prove 
Theorem \ref{GorC}. In Section \ref{section7}, we introduce the notation and review some facts concerning Kähler differential modules and differential modules over rings and schemes. Subsequently, in Section \ref{section8}, we explore the properties of the Kähler differential module $\Omega_{X/Y}^{(n)}$, focusing on its reflexivity, $k$-torsionfreeness, and $k$-syzygy properties, particularly for $k = 1, 2$. Finally, in the same section, we offer a partial answer to the generalized Herzog-Vasconcelos conjecture (GHVC).

\section{Setup and preliminaries}\label{section2}
Throughout this paper,  $R$ is a commutative Noetherian ring, and all $R$-modules are considered to be finitely generated. 
\begin{definition}\label{definitions} Let $M$ and $N$ be $R$-modules, and let $k$ be a non-negative integer.
\begin{enumerate}
\item We set $M^\ast=\Hom_R(M,R)$ and $M^{\ast \ast}=(M^\ast)^\ast$.
    \item Let $$\boldsymbol{P}:\cdots \longrightarrow P_i \stackrel{\varphi_i}{\longrightarrow} P_{i-1} \longrightarrow \cdots \longrightarrow P_1 \stackrel{\varphi_1}{\longrightarrow} P_0 \stackrel{\varphi_0}{\longrightarrow} M \longrightarrow 0$$ be a projective resolution of $M$.  
    \begin{enumerate}
        \item[(2.a)]  For $k\geq 1$, the $k$-\textit{syzygy} of $M$, denoted by $\Omega^k(M)$, is defined as the kernel of the map $\varphi_{k-1}$. When $k=0$, we set $\Omega^0(M)=M$.
        \item[(2.b)] The \textit{Auslander transpose} of $M$, denoted by $\operatorname{Tr}(M)$, is defined as the cokernel of  the induced map $\varphi_{1}^\ast: P_0^\ast \to P_1^\ast.$ When $k \geq 1$, we set  $$\mathcal{T}_k(M)=\operatorname{Tr}(\Omega^{k-1}(M)).$$
    \end{enumerate}
    \item We write $M \approx N$ if there exist projective $R$-modules $F,G$ such that $M \oplus F \cong N \oplus G$. 
   \item For $k\geq 1$, $M$ is said to be a $k$-\textit{syzygy} if there exists an exact sequence $$0 \rightarrow M \rightarrow P_{k-1} \rightarrow P_{k-2} \rightarrow \cdots \rightarrow P_1 \rightarrow P_0,$$ where each $P_i$ is  projective. For $k=0$, every $R$-module is said to be a $0$-\textit{syzygy}.  
   \item $M$ is $k$-\textit{torsionfree} if $\Ext_R^i(\operatorname{Tr}(M),R)=0$ for all $1\leq i \leq k$. 
   \item We set 
$X^k(R)$ (resp. $\widetilde{X}^k(R)$) to the set of all prime  ideals $\mathfrak{p}$ of $R$   such that $\operatorname{ht}(\mathfrak{p})\leq k$ (resp. $\operatorname{depth}_{R_\mathfrak{p}} (R_\mathfrak{p})\leq k)$.
\item We say that $M$ satisfies $(S_k)$ (resp. $(\widetilde{S}_k)$) if $\operatorname{depth}_{R_\mathfrak{p}}(M_\mathfrak{p})\geq \min\{k, \operatorname{ht}(\mathfrak{p})\}$ (resp. $\operatorname{depth}_{R_\mathfrak{p}}(M_\mathfrak{p})\geq \min\{k, \depth_{R_\mathfrak{p}}(R_\mathfrak{p})\} )$ for all prime ideals $\mathfrak{p}$ of $R$. 
\item We say that  $R$ satisfies $(G_k)$ if the local ring $R_\mathfrak{p}$ is Gorenstein for all $\mathfrak{p} \in X^k(R)$. 
   
\end{enumerate}
\end{definition}

The $R$-modules $\operatorname{Tr}(M)$ and $\Omega^k(M)$ ($k\geq 0$) are uniquely determined up to projective summands. It is easy to see that $\operatorname{Tr}(\operatorname{Tr}(M))\approx M$. Note that $M$ is a $k$-syzygy if and only if there exists an $R$-module $L$ such that $M \approx \Omega^k(L)$. Whenever $R$ is local, we consider the Auslander transpose and the syzygies (of an $R$-module) defined using minimal free resolutions. In this case, these $R$-modules are defined uniquely up to isomorphism (rather than projective equivalence).

The notion of Gorenstein dimension was introduced by Auslander \cite{AnneauxdeGorensteinettorsionenalgebrecommutative}  and developed by Auslander and Bridger in \cite{Stablemoduletheory}. 
\begin{definition} Let $M$ be an $R$-module. 
 \begin{enumerate}
     \item We say that $M$ is  \textit{totally reflexive} if  $M$ is reflexive and $\Ext_R^i(M,R)=\Ext_R^i(M^\ast,R)=0$ for all $i>0$. 
     \item The \textit{Gorenstein dimension} of $M$, denoted by $\operatorname{G-dim}_R(M)$, is defined to be the infimum of all non-negative integers $k$ such that there exists an exact sequence  
$$
0 \rightarrow G_k \rightarrow \cdots \rightarrow G_0 \rightarrow M \rightarrow 0,
$$
where each $G_i$ is totally reflexive. 
 \end{enumerate}
\end{definition}
We can observe that $\operatorname{G-dim}_R(M)=0$ if and only if $M$ is totally reflexive. Below, we collect some facts related to the Gorenstein dimension.

\begin{facts}\label{eqpl} Let $M$ be an $R$-module.
\begin{enumerate}
    \item (\cite[Theorem 1.2.7]{GorensteinDimensions}) If $M\not=0$ and $\operatorname{G-dim}_R(M)<\infty$, then $$\operatorname{G-dim}_R(M)=\sup\{i\geq 0: \Ext_R^i(M,R)\not=0\}.$$
    \item (\cite[Lemma 3.19(1)]{Stablemoduletheory}) $M$ is totally reflexive if and only if   $\operatorname{Tr} (M)$ is totally reflexive. 
    \item (\cite[Corollary 1.2.9]{GorensteinDimensions}) Let $0 \to M \to N \to L \to 0$ be an exact sequence of $R$-modules. If two $R$-modules of the sequence have finite Gorenstein dimension, then so has the third. 
    \item (\cite[Proposition 1.2.10]{GorensteinDimensions}) $\operatorname{G-dim}_R(M)\leq \operatorname{pd}_R(M)$ with equality if $\operatorname{pd}_R(M)<\infty$. 
\end{enumerate}
Now, assume that $R$ is local. 
\begin{enumerate}
    \item[(5)] (\cite[Theorem 1.4.8]{GorensteinDimensions}) If  $M\not=0$ and  $\operatorname{G-dim}_R(M)<\infty$, then    $$\operatorname{G-dim}_R(M)+ \operatorname{depth}_R(M)=\operatorname{depth}_R(R).$$    In particular, $\operatorname{G-dim}_R(M)\leq \operatorname{depth}_R(R)$.
    \item[(6)] (\cite[Theorem 1.4.9]{GorensteinDimensions}) Let $k$ be the residual field of $R$. The following conditions are equivalent:
     \begin{enumerate}[(6.1)]
         \item $R$ is Gorenstein.
         \item $\operatorname{G-dim}_R(M)<\infty$ for all (finitely generated) $R$-modules $M$.
         \item $\operatorname{G-dim}_R(k)<\infty$.
     \end{enumerate}
\end{enumerate}
\end{facts}
The formula given in Facts \ref{eqpl}(5) is known as the \textit{Auslander-Bridger formula}. 
From the definition of the Auslander transpose and Facts \ref{eqpl}(3) we can derive the following remark.
\begin{remark}\label{wss}
    Let $R$ be a ring and let $M$ be an $R$-module. Then $\operatorname{G-dim}_R(M^\ast)<\infty$ if and only if $\operatorname{G-dim}_R(\operatorname{Tr}(M))<\infty$.
\end{remark}

\section{A generalization of a theorem of Auslander and Bridger}
The main goal of this section is to provide a generalization of a celebrated result that Evans and Griffith attribute  Auslander and Bridger in \cite[Theorem 3.8]{Syzygies} (see also \cite[Lemma 1.3]{AscentOfFiniteCohen–MacaulayType}), which characterizes a module to be $k$-torsionfree, $k$-syzygy, involving the condition $(\widetilde{S}_k)$.

Before presenting such a generalization, we provide the following result, which will be crucial for the proof of our main theorem in this section.

\begin{lemma}\label{syad01}
Let $R$ be a local ring of depth $t$, $M$ be a non-zero $R$-module and $0\leq  k \leq t$ be an integer.  If $M$ is a $k$-syzygy of an $R$-module of depth $0$, then $\operatorname{depth}_R(M)=k$.
\end{lemma}
\begin{proof}
    It is sufficient to consider $k \geq 1$. We proceed by induction on $k$. Suppose $k=1$. Then $t\geq 1$ and there exists an exact sequence
    $0 \to M \to F_0 \to L\to 0,$
    where $F_0$  is free and $L$ has depth $0$. The depth lemma gives us the inequalities
    \begin{equation}\label{deax}
         0=\operatorname{depth}_R(L)\geq \min\{\operatorname{depth}_R(M)-1, t \}
    \end{equation}
    and  \begin{equation}\label{deax2}
        \operatorname{depth}_R(M)\geq \min\{t,1\}=1. 
     \end{equation}
         Since $t\geq 1$,  it follows from  \eqref{deax} that $\operatorname{depth}_R(M)\leq 1$. Thus  \eqref{deax2} yields that $\operatorname{depth}_R(M)=1$.

    Now suppose $1<k \leq t$. From the given hypotheses, there exists an exact sequence 
    $$0 \to M \to F_{k-1} \to \cdots \to F_{0} \to L \to 0,$$
    where each $F_i$ is free and $\operatorname{depth}_R(L)=0$. Then it induces exact sequences
\begin{equation}\label{poa1}
    0 \to M \to F_{k-1} \to N \to 0
\end{equation}
    and 
    \begin{equation}\label{poa2}
        0 \to N \to F_{k-2} \to \cdots \to F_0 \to L \to 0.
    \end{equation}

    If $N=0$, then \eqref{poa2} says that  $\operatorname{pd}_R(L)\leq k-2$ and, by the Auslander-Buchsbaum formula, $t\leq k-2,$ which does not occurs by assumption. So, $N\not=0$. Then by induction, we have that $\operatorname{depth}_R(N)=k-1$. Applying the depth lemma in \eqref{poa1} we get that 
    
    $$\operatorname{depth}_R(M)\geq \min\{t, \operatorname{depth}_R(N)+1\}= \min \{t, k\} =k$$
    and  
    $$k-1=\operatorname{depth}_R(N)\geq \min\{ \operatorname{depth}_R(M)-1, t\}.$$
Thus, if $t\leq \operatorname{depth}_R(M)-1$, then $k-1\geq t$, which contradicts our assumption. Therefore, we must have $t > \operatorname{depth}_R(M)-1$. Consequently,  $\operatorname{depth}_R(M)-1\leq k-1$, which implies $\operatorname{depth}_R(M)\leq k$. Therefore, $k=\operatorname{depth}_R(M)$.
       \end{proof}
Below, we provide a characterization of a module to be $k$-torsionfree or $k$-syzygy, involving the condition $(\widetilde{S}_k)$. This characterization generalizes \cite[Theorem 3.8]{Syzygies}

\begin{theorem}\label{GteosPa} Let $R$ be a ring  and let $k$ be a non-negative integer. Let $M$ be an $R$-module such that $M^\ast$ has locally finite Gorenstein dimension on $\widetilde{X}^{k-1}(R)$. Then the following conditions are equivalent:
    \begin{enumerate}
        \item $M$ is $k$-torsionfree.
        \item $M$ is $k$-syzygy.
        \item $M$ satisfies $(\widetilde{S}_k)$. 
    \end{enumerate}
\end{theorem}
\begin{proof}
The implications $(1) \Rightarrow (2) \Rightarrow (3)$ follow from \cite[Proposition 38]{GorensteinDimensionAndTorsionOfModulesOverCommutativeNoetherianRings}. Now let us prove $(3)\Rightarrow (1)$. We prove this by induction on $k$. If $k=0$, there is nothing to prove. Assume that $k\geq 1$. By induction, $M$ is $(k-1)$-torsionfree, that is  $\Ext_{R}^{i}(\operatorname{Tr}(M),R)=0$ for all $1\leq i \leq k-1$, and so it remains to conclude that $\Ext_R^{k}(\operatorname{Tr}(M),R)=0$. 

 By contradiction, assume that $\Ext_R^k(\operatorname{Tr}(M), R)\not=0$. Let $\mathfrak{p} \in \operatorname{Ass}_{R}(\Ext_R^k(\operatorname{Tr}(M), R))$.  Then 
$
\mathfrak{p} R_{\mathfrak{p}} \in \operatorname{Ass}_{R_{\mathfrak{p}}}\left(\operatorname{Ext}_{R_{\mathfrak{p}}}^k\left(\operatorname{Tr}_{R_{\mathfrak{p}}} (M_{\mathfrak{p}}), R_{\mathfrak{p}}\right)\right)
$. For a Noetherian local ring $S$ and a (finitely generated) $S$-module $N$,  the $S$-module $\Ext_S^1(N,S)$ is a submodule of $\operatorname{Tr}(N)$. Thus, taking $S=R_\mathfrak{p}$ and  $N=\Omega_{R_\mathfrak{p}}^{k-1}(\operatorname{Tr}_{R_\mathfrak{p}}(M_\mathfrak{p}))$, we obtain an exact sequence 
    $$0 \to \Ext_{R_\mathfrak{p}}^k\left(\operatorname{Tr_{R_\mathfrak{p}}}(M_\mathfrak{p}),R_\mathfrak{p}\right) \to \mathcal{T}_k(\operatorname{Tr}_{R_\mathfrak{p}}(M_\mathfrak{p}) ), $$
    which implies that $\mathfrak{p}R_\mathfrak{p} \in \operatorname{Ass}_{R_\mathfrak{p}}(\mathcal{T}_k(\operatorname{Tr}_{R_\mathfrak{p}}(M_\mathfrak{p}))$, and hence $$\operatorname{depth}_{R_\mathfrak{p}}(\mathcal{T}_k(\operatorname{Tr}_{R_\mathfrak{p}}(M_\mathfrak{p}))=0.$$ Since $\Ext_{R_\mathfrak{p}}^i(\operatorname{Tr}_{R_\mathfrak{p}}(M_\mathfrak{p}), R_\mathfrak{p})=0$ for all $1\leq i \leq k-1$, we see from \cite[5.6]{HomologicalDimensionsOfRigidModules} that $$M_\mathfrak{p} \approx \Omega_{R_\mathfrak{p}}^{k-1}\mathcal{T}_k(\operatorname{Tr}_{R_\mathfrak{p}}(M_\mathfrak{p})).$$

    We claim that $k\leq \depth_{R_\mathfrak{p}}(R_\mathfrak{p})$. Indeed, if $k>\operatorname{depth}_{R_\mathfrak{p}}(R_\mathfrak{p})$, then by assumption $\operatorname{G-dim}_{R_\mathfrak{p}}(\operatorname{Tr}_{R_\mathfrak{p}}(M_\mathfrak{p}))<\infty$, and by the Auslander-Bridger formula, $\operatorname{G-dim}_{R_\mathfrak{p}}(\operatorname{Tr}_{R_\mathfrak{p}}  (M_\mathfrak{p}))<k$. But since $\Ext_{R_\mathfrak{p}}^i(\operatorname{Tr}_{R_\mathfrak{p}}(M_\mathfrak{p}), R_\mathfrak{p})=0$ for all $1\leq i \leq k-1$, we get that $\operatorname{Tr}_{R_\mathfrak{p}}(M_\mathfrak{p})$ is a totally reflexive. This shows that $\Ext_{R_\mathfrak{p}}^k(\operatorname{Tr}_{R_\mathfrak{p}}(M_\mathfrak{p}),R_\mathfrak{p})=0$, which does not occurs.  

    Now, as $M_\mathfrak{p} \approx \Omega_{R_\mathfrak{p}}^{k-1}\mathcal{T}_k(\operatorname{Tr}_{R_\mathfrak{p}}(M_\mathfrak{p})) $, it follows from Lemma \ref{syad01} that $\operatorname{depth}_{R_\mathfrak{p}}(M_\mathfrak{p})=k-1$. But since $M$  satisfies $(\widetilde{S}_k)$, we get that $\operatorname{depth}_{R_\mathfrak{p}}(M_\mathfrak{p})\geq \min\{k,\operatorname{depth}_{R_\mathfrak{p}}(R_\mathfrak{p})\}\linebreak=k$, a contradiction. 
\end{proof}
Observe that Theorem \ref{GteosPa} is the version with the dual of an $R$-module being of finite Gorenstein dimension of \cite[Proposition 2.4]{LinkageOfModulesAndTheSerreConditions} when $C=R$. Now, we explore some consequences of Theorem \ref{GteosPa}.
\begin{corollary}\label{saq54}
    Let $R$ be a ring and let $M$ be an $R$-module. Suppose that $\operatorname{G-dim}_{R_\mathfrak{p}}(M_\mathfrak{p}^\ast)<\infty$ for all $\mathfrak{p} \in \operatorname{Ass}_R(R)$. Then the following conditions are equivalent:
    \begin{enumerate}
        \item $M$ is $1$-torsionfree.
        \item $M$ is $1$-syzygy.
        \item $M$ satisfies $(\widetilde{S}_1)$.
    \end{enumerate}
\end{corollary}

\begin{corollary}\label{saq55} 
    Let $R$ be a ring and let $M$ be an $R$-module. Suppose that $\operatorname{G-dim}_{R_\mathfrak{p}}(M_\mathfrak{p}^\ast)<\infty$ for all $\mathfrak{p} \in \widetilde{X}^1(R)$. Then the following conditions are equivalent:
    \begin{enumerate}
        \item $M$ is $2$-torsionfree.
        \item $M$ is $2$-syzygy.
        \item $M$ satisfies $(\widetilde{S}_2)$.
        \item $M$ is reflexive.
        \end{enumerate}
\end{corollary}

As another consequence, this theorem allows us to recover \cite[Theorem 5.8(1)]{HomologicalDimensionsOfRigidModules}. Additionally, it provides an answer to the following question posed in \cite{WhenArenSysymodulesnTorsionfree}. 

\begin{question}[{\cite[Question 1.1]{WhenArenSysymodulesnTorsionfree}}]
    When $k$-syzygy modules are $k$-torsionfree?
\end{question}

The following result is an analogous version of \cite[Theorem 4.5]{ktorsionlessmoduleswithfiniteGorensteindimension}, which is one of the main results of that paper.
\begin{theorem}
    \label{aqps} Let $R$ be a ring, and let $M$ be an $R$-module such that $M^\ast$ has locally finite Gorenstein dimension on $\widetilde{X}^{k-1}(R)$.  Then $M$ is $k$-torsionfree if and only if the following hold:
    \begin{enumerate}[(1)]
        \item $M_\mathfrak{p}$ is $(k-1)$-torsionfree for all $\mathfrak{p} \in \widetilde{X}^{k-1}(R)$. 
        \item $\operatorname{depth}_{R_\mathfrak{p}}(M_\mathfrak{p})\geq k$ for all $\mathfrak{p} \in \Spec(R)$ with $\operatorname{depth}_{R_\mathfrak{p}}(R_\mathfrak{p})\geq k$. 
    \end{enumerate}
\end{theorem}
\begin{proof}
    If $M$ is $k$-torsionfree, then it is clear that (1) is true. Besides, by Theorem \ref{GteosPa}, $M$ satisfies $(\widetilde{S}_k)$, which shows that (2) holds.

    Reciprocally, suppose that (1) and (2) hold. By Theorem \ref{GteosPa}, it is sufficient to show that $M$ satisfies $(\widetilde{S}_k),$ that is $\operatorname{depth}_{R_\mathfrak{p}}(M_\mathfrak{p})\geq \min\{k, \operatorname{depth}_{R_\mathfrak{p}}(R_\mathfrak{p})\}$ for all $\mathfrak{p} \in \Spec(R)$. Let $\mathfrak{p} \in \Supp (M)$.  By (2), it holds if $\operatorname{depth}_{R_\mathfrak{p}}(R_\mathfrak{p})\geq k$. Now, let $\depth_{R_\mathfrak{p}}(R_\mathfrak{p})\leq k-1$. Then by the Auslander-Bridger formula, $\operatorname{G-dim}_{R_\mathfrak{p}}(\operatorname{Tr}_{R_\mathfrak{p}}(M_\mathfrak{p}))\leq k-1$. Thus, by (1),  $\operatorname{Tr}_{R_\mathfrak{p}}(M_\mathfrak{p})$ is totally reflexive, and hence $M_\mathfrak{p}$ as well.  In particular, $\operatorname{depth}_{R_\mathfrak{p}}(M_\mathfrak{p})=\operatorname{depth}_{R_\mathfrak{p}}(R_\mathfrak{p})=\min\{k, \operatorname{depth}_{R_\mathfrak{p}}(R_\mathfrak{p})\}$. 
\end{proof}

\section{Characterization of totally reflexives}
In this section, $R$ is a local ring. Next,  we derive some criteria for an $R$-module $M$ to be  totally reflexive, assuming the condition $\operatorname{G-dim}_R(M^\ast)<\infty$ and involving the vanishing of the Ext functor. We first provide a general criterion via the vanishing of $\Ext_R^{1\leq i \leq \operatorname{depth}_R(R)}(M,R)$, and subsequently, we provide criteria over Cohen-Macaulay local rings and later for locally totally reflexive  modules on the punctured spectrum. 

\begin{theorem}  \label{pos5} Let $R$ be  a local ring of depth $t$. Let $M$ be an $R$-module such that $\operatorname{G-dim}_R(M^\ast)<\infty$. Then the following are equivalent:
\begin{enumerate}
    \item $M$ satisfies $(\widetilde{S}_t)$.
    \item $M$ is totally reflexive.
    \item $\Ext_R^i(M,R)=0$ for all $1\leq i \leq t$. 
\end{enumerate}
\end{theorem}
\begin{proof}
    $(1) \Rightarrow (2)$. According to Theorem \ref{GteosPa} or \cite[Theorem 5.8(1)]{HomologicalDimensionsOfRigidModules}, $M$ is $t$-torsionfree, meaning that $\Ext_R^i(\operatorname{Tr}(M),R)=0$ for all $1\leq i \leq t$. Since $\operatorname{G-dim}_R(\operatorname{Tr}(M))\leq t,$ it follows that $\operatorname{Tr}(M)$ is totally reflexive. Hence, $M$ is totally reflexive. 
    
    $(2) \Rightarrow (3)$. This follows from the definition of totally reflexive modules.
     
     $(3) \Rightarrow (1)$.  Suppose that $\Ext_R^i(M,R)=0$ for all $1\leq i \leq t$. Since $M \approx \operatorname{Tr}(\operatorname{Tr}(M))$, this implies that $\operatorname{Tr}(M)$ is $t$-torsionfree. Given that $\operatorname{G-dim}_R(\operatorname{Tr}(M))<\infty$, according to \cite[Theorem 42]{GorensteinDimensionAndTorsionOfModulesOverCommutativeNoetherianRings}, $\operatorname{Tr}(M)$ satisfies $(\widetilde{S}_t)$. Thus, $\operatorname{depth}_R(\operatorname{Tr}(M))\geq t$ since $\depth_R(R)=t$. By the Auslander-Bridger formula, $\operatorname{Tr}(M)$ is totally reflexive. Thus, \linebreak $\Ext_R^i(\operatorname{Tr}(M),R)=0$ for all $i>0$. In particular, $M$ is $t$-torsionfree, and by \cite[Proposition 11]{GorensteinDimensionAndTorsionOfModulesOverCommutativeNoetherianRings}, $M$ satisfies $(\widetilde{S}_t)$.   
\end{proof}

\subsection*{Over Cohen-Macaulay local rings} As an immediate application of Theorem \ref{pos5}, we obtain the following characterization for a module over a Cohen-Macaulay local ring to be totally reflexive, involving the vanishing of the Ext functor. 
\begin{corollary}
    \label{Mc5s}
     Let $R$ be a Cohen-Macaulay local ring of dimension $d$, and let $M$ be a non-zero $R$-module such that $\operatorname{G-dim}_R(M^\ast)<\infty$. Then the following are equivalent:
     \begin{enumerate}
         \item $M$ is maximal Cohen-Macaulay.
         \item $M$ is totally reflexive.
         \item $\Ext_R^i(M,R)=0$ for all $1\leq i \leq d$.
     \end{enumerate}
\end{corollary}

\begin{remark} It is known that if $R$ is a Gorenstein local ring and $M$ is maximal Cohen-Macaulay, then $\Ext_R^i(M,R)=0$ for all $i>0$. We can observe from Corollary \ref{Mc5s}, that in this fact instead of $R$ being Gorenstein, we can consider the weaker condition of $\operatorname{G-dim}_R(M^\ast)<\infty$. 
\end{remark}

We recall that the \textit{grade} of an $R$-module $M$ is defined as: $$\operatorname{grade}_R(M)=\inf\{i\geq 0: \Ext_R^i(M,R)\not=0\}. $$
\begin{corollary}
    \label{sdq}
    Let $R$ be a Cohen-Macaulay local ring, and let $M$ be a Cohen-Macaulay $R$-module such that $M^\ast \not=0$ and $\operatorname{G-dim}_R(M^\ast)<\infty$. Then $M$ is totally reflexive.  
\end{corollary}
\begin{proof}
 Since $M^\ast\not=0$, note that $\operatorname{grade}_R(M)=0$. Moreover, by \cite[Corollary 2.1.4]{bruns}, $\operatorname{grade}_R(M)=\operatorname{dim}(R)-\operatorname{dim}_R(M)$ and hence $\operatorname{dim}_R(M)=\operatorname{dim}(R).$ Thus from Cohen-Macaulayness of $M$, we see  that $M$ is maximal Cohen-Macaulay. Therefore, Corollary \ref{Mc5s} asserts that $M$ is totally reflexive.  
\end{proof}

Now, as an application of  Corollary \ref{Mc5s}, we  extend the formula presented in \cite[Lemma 4.1]{LinkageOfFiniteGorensteinDimensionModules} to all $R$-modules $M$ such that $\operatorname{G-dim}_R(M^\ast)<\infty$.  But before, we recall the definition of reduced grade. The \textit{reduced grade} of an $R$-module $M$ is defined as: 
$$\operatorname{r.grade}_R(M)=\inf\{i>0: \Ext_R^i(M,R)\not=0\}.$$

\begin{remark} \label{grade} Let $R$ be a ring and let $M$ be an $R$-module. If $\operatorname{grade}_R(M)>0,$ then $\operatorname{grade}_R(M)=\operatorname{r.grade}_R(M).$
\end{remark}

\begin{proposition} \label{sepxw}
    Let $R$ be a Cohen-Macaulay local ring of dimension $d$, and let $M$ be an $R$-module such that $\operatorname{G-dim}_R(M^\ast)<\infty$. Then the equality
    \begin{equation}\label{ewbn}
    \sup \{i \geq 0 :H_{\mathfrak{m}}^i(M) \neq 0 \text{ and } i \neq d\}=d-\operatorname{r.grade}(M)
\end{equation}
holds in $\mathbb{Z} \cup \{\pm \infty\}$.
\end{proposition}
\begin{proof}
     We may assume that $R$ is complete, and hence $R$ has a canonical module $\omega$.  We prove the equality by considering the following cases:
     \begin{enumerate}
         \item $M$ is not maximal Cohen-Macaulay and $\operatorname{dim}_R(M)=\operatorname{dim}(R)$.
        \item $M$ is not maximal Cohen-Macaulay and $\operatorname{dim}_R(M)\not=\operatorname{dim}(R)$. 
        \item $M$ is maximal Cohen-Macaulay. 
     \end{enumerate}
     The equality in the case (1) was proved in \cite[Lemma 4.1]{LinkageOfFiniteGorensteinDimensionModules}. 
     
    Now, consider the case (2). Since $\operatorname{dim}_R(M)\not=\dim(R)$, we see from \cite[Theorem 3.5.7]{bruns} that $\sup\{i \geq 0: H_\mathfrak{m}^i(M)\not=0 \text{ and } i \not=d\}=\operatorname{dim}_R(M).$  Moreover, by \cite[Corollary 2.1.4]{bruns}, we get  that  $\operatorname{grade}_R(M)=d-\operatorname{dim}_R(M)>0$. Then by Remark \ref{grade}, $\operatorname{r.grade}_R(M)=d-\operatorname{dim}_R(M)$. Now observe that the desired equality holds. 

     Finally, let us prove the equality in the case (3). In this case, by \cite[Theorem 3.5.7]{bruns}, $\sup\{i \geq 0: H_\mathfrak{m}^i(M)\not=0 \text{ and }  i \not=d  \}=-\infty$. On the other hand, by Corollary \ref{Mc5s}, $M$ is totally reflexive and hence  $\Ext_R^i(M,R)=0$ for all $i\geq 1$. Thus $\operatorname{r.grade}_R(M)=\infty$ and we can see that the required equality holds.  
    \end{proof}  
     
\begin{corollary}\label{sepx}
     Let $R$ be a Cohen-Macaulay local ring with a canonical module $\omega$, and let $M$ be an $R$-module such that $\operatorname{G-dim}_R(M^\ast)<\infty$. Then: $$\operatorname{r.grade}(M)=\inf\{i >0: \Ext_R^i(M,\omega)\not=0\}.$$
\end{corollary}
\begin{proof}
We may assume that $R$ is complete. Because of \cite[Theorem 3.5.8]{bruns}, the desired equality is a reformulation of the equality \eqref{ewbn}.  
\end{proof}

\begin{theorem}\label{genCM}
    Let $R$ be a Cohen-Macaulay local ring of dimension $d$, and let $M$ be an $R$-module such that $\operatorname{G-dim}_R(M^\ast)<\infty$. Let $n$ be a non-negative integer such that $\operatorname{depth}_R(M)\geq n$. If $\Ext_R^i(M,R)=0$ for all $1\leq i \leq d-n,$ then $M$ is totally reflexive.  
\end{theorem}
\begin{proof}
We may assume that $M\not=0$ and (passing to the completion of $R$ if necessary) that $R$ has a canonical module $\omega$. 

Since $\Ext_R^i(M,R)=0$ for all $1\leq i \leq d-n$ and $\operatorname{G-dim}_R(M^\ast)<\infty$, by Corollary \ref{sepx}, $\Ext_R^i(M,\omega)=0$ for all $1\leq i \leq d-n$. Since  \cite[Exercise 3.1.24]{bruns} says us that $$d-\operatorname{depth}_R(M)=\sup\{ i\geq 0: \Ext_R^i(M,\omega)\not=0\},$$
then  $d-\operatorname{depth}_R(M)=0$ or $d-\operatorname{depth}_R(M)>d-n$. The second case does not occur because $\operatorname{depth}_R(M)\geq n$ by assumption. Thus $d=\operatorname{depth}_R(M)$ and hence $M$ is maximal Cohen-Macaulay. Thus, by Corollary \ref{Mc5s},  $M$ is totally reflexive. 

\end{proof}

\subsection*{Over locally totally reflexive modules on the punctured spectrum} The following result provides a totally reflexivity criterion for $R$-modules $M$ with $\operatorname{G-dim}_R(M^\ast)<\infty$  that are locally totally reflexives on the punctured spectrum of $R$. 
 \begin{theorem}\label{seaa1}
     Let $R$ be a local ring of depth $t$, and let $M$ be a non-zero $R$-module such that $\operatorname{G-dim}_R(M^\ast)<\infty$. Suppose that $M$  is locally totally reflexive on the punctured spectrum of $R$. Then $\depth_R(M)\leq t$ with equality if and only if $M$ is totally reflexive.
 \end{theorem}
 \begin{proof}
     First, we prove that $\operatorname{depth}_R(M)\leq t$. Set $r=\operatorname{depth}_R(M)$ and suppose that $r>t$. Then $M$ satisfies $(\widetilde{S}_r)$ since $M$ is locally totally reflexive on the punctured spectrum of $R$. In particular, $M$ satisfies $(\widetilde{S}_t)$. By Theorem \ref{pos5}, $M$ is totally reflexive, which contradicts the fact that $r>t$.

     Now, we prove the second part of the result. It is clear that if $M$ is totally reflexive, then $\operatorname{depth}_R(M)=t$. Suppose now that $\operatorname{depth}_R(M)=t$. Since $M$ is locally totally reflexive on the punctured spectrum of $R$ and $\operatorname{depth}_R(M)=t$, note that $M$ satisfies $(\widetilde{S}_t)$. Hence, by Theorem \ref{pos5}, $M$ is totally reflexive.
 \end{proof}

\begin{corollary}\label{a11p}
    Let $R$ be a local ring of depth $t$, and let $M$ be a non-zero $R$ module. Suppose that   $M$ is locally totally reflexive on the punctured spectrum of $R$. Suppose that $\operatorname{G-dim}_R(M^\ast)<\infty$. Then the following are equivalent:   \begin{enumerate}
    \item $\operatorname{depth}_R(M)\geq t$.
        \item $\operatorname{depth}_R(M)=t$.
        \item $M$ is totally reflexive.
        \item $\Ext_R^i(M,R)=0$ for all $1\leq i \leq t$.
    \end{enumerate}
\end{corollary}
\begin{proof}
This follows from Theorems \ref{seaa1} and \ref{pos5}.
\end{proof}

In this section, we have derived a number of criteria for an $R$-module $M$  with $\operatorname{G-dim}_R(M^\ast)<\infty$ to be totally reflexive. Thus, we close this section with the following question. 
\begin{question}\label{sqeeww}
    Let $R$ be a ring and let $M$ be an $R$-module. When does $\operatorname{G-dim}_R(M^\ast)<\infty$  imply that $\operatorname{G-dim}_R(M)<\infty$? 
\end{question}

\section{Freeness Criteria}

In this section, we aim to obtain freeness criteria by developing the following topics:

\begin{enumerate}[(I)]
    \item To derive freeness criteria from Theorem \ref{genCM}.
    \item To provide generalizations of a celebrated result of Araya (Theorem \ref{araq}) with the condition that the dual of a module is of finite Gorenstein dimension.
    \item To discuss the Auslander-Reiten conjecture for $R$-modules $M$ such that \break $\operatorname{G-dim}_R(M^\ast)<\infty$ and $\operatorname{pd}_R(\Hom_R(M,M))<\infty$.  
\end{enumerate}

\subsection*{Freeness criteria from Theorem \ref{genCM}} As an application of Theorem \ref{genCM}, we obtain the following result: 

\begin{theorem}\label{pesge}
     Let $R$ be a Cohen-Macaulay local ring of dimension $d$, and let $M$ be an $R$-module such that $\operatorname{pd}_R(M^\ast)<\infty$. Let $n$ be a non-negative integer such that $\operatorname{depth}_R(M)\geq n$ and $\Ext_R^i(M, R)=0$ for all $1\leq i \leq d-n$. Then $M$ is free.
 \end{theorem}
 \begin{proof}
    As $\operatorname{pd}_R(M^\ast)<\infty$, then $\operatorname{G-dim}_R(M^\ast)<\infty$. By Theorem \ref{genCM}, $M$ is totally reflexive. Thus $M$ is reflexive and $\operatorname{pd}_R(M^\ast)=\operatorname{G-dim}_R(M^\ast)=0$. Therefore, $M$ is free.
  \end{proof}
  \begin{corollary} \label{spkp}
      Let $R$ be a Cohen-Macaulay local ring, and let $M$ be an $R$-module such that $\operatorname{pd}_R(M^\ast)<\infty$. If $M$ is maximal Cohen-Macaulay, then $M$ is free.
  \end{corollary}
  Similarly to the proof of Corollary \ref{sdq}, we obtain the following. 
  \begin{corollary}\label{erwq}
      Let $R$ be a Cohen-Macaulay local ring, and let $M$ be an $R$-module such that $M^\ast \not=0$ and  $\operatorname{pd}_R(M^\ast)<\infty$. If $M$ is Cohen-Macaulay, then $M$ is free.
  \end{corollary}

As an application of Corollary \ref{spkp}, we will provide an answer to the Generalization of Herzog-Vasconcelos's conjecture (see Proposition \ref{diff1}). Another application is given below. 
 \begin{theorem}
     Let $(R,\mathfrak{m},k)$ be a Cohen-Macaulay local ring. The maximal ideal  $\mathfrak{m}$ is Cohen-Macaulay if and only $\operatorname{dim}(R)\leq 1$. 
 \end{theorem}
 \begin{proof} Assume that $\mathfrak{m}$ is Cohen-Macaulay and that $\operatorname{dim}(R)=\operatorname{depth}_R(R)\geq 2$. Consider the exact sequence 
\begin{equation}\label{setx}
    0 \to \mathfrak{m} \to R \to k \to 0.
\end{equation}
It induces an exact sequence 
$$0 \to k^\ast \to R \to \mathfrak{m}^\ast \to \Ext_R^1(k,R).$$
Since  $\operatorname{depth}_R(R)\geq 2,$ then $k^\ast=0=\Ext_R^1(k,R)$ and $\mathfrak{m}^\ast \cong R$.  Therefore, by Corollary \ref{erwq},  $\mathfrak{m}$ is free and hence $\mathfrak{m} \cong R$. From the exact sequence \eqref{setx}, we obtain an exact sequence 
$$0 \to \Hom_R(k, \mathfrak{m}) \to \Hom_R(k,R) \to \Hom_R(k,k) \to \Ext_R^1(k,\mathfrak{m}).$$ 
 As $\operatorname{depth}_R(R)\geq 2$ and $\mathfrak{m}\cong R$,  then $\Hom_R(k,\mathfrak{m}) \cong \Hom_R(k, R)=0=\Ext_R^1(k,R)\cong \Ext_R^1(k,\mathfrak{m})$. Consequently, $\Hom_R(k,k)=0,$ which does not occurs.
 
 Conversely, assume that $\operatorname{dim}(R)\leq 1$. If $\dim(R)=0$ it is clear that all $R$-modules are  Cohen-Macaulay. Suppose $\operatorname{dim}(R)=1$. Since $R$ is Cohen-Macaulay, then there exists an $R$-regular element $x \in \mathfrak{m}$. Note that $x$ is also an $\mathfrak{m}$-regular element. Therefore $1\leq \operatorname{depth}_R(\mathfrak{m})\leq \operatorname{dim} (R)=1,$ which shows that $\mathfrak{m}$ is maximal Cohen-Macaulay.
 \end{proof}

\subsection*{Some generalizations of a result of Araya} As mentioned in the introduction, after the publication of the Araya's result referenced in Theorem \ref{araq}, results have appeared that generalize or recover such theorem, first in the context of Gorenstein rings and later in more general contexts. Motivated by this, we provide some generalizations of that celebrated result in terms of the dual of an $R$-module is of finite Gorenstein dimension.
 
\begin{proposition}\label{r61}
    Let $R$ be a local ring of depth $t\geq 1$. Let $M$ be a locally free $R$-module  on the punctured spectrum such that $\operatorname{G-dim}_R(M^\ast)<\infty$. If $\Ext_R^{t-1}(M,M)=0$ and  $\Ext_R^i(M,R)=0$ for all $1\leq i \leq t$, then $M$ is free. 
\end{proposition}
\begin{proof}
If $t=1$, then $\Hom_R(M,M)\cong \Ext_R^0(M,M)=0$, whence $M$ is zero and hence free. Suppose that $t\geq 2$. By Theorem \ref{pos5}, $M$ is totally reflexive and hence $\operatorname{Tr}(M)$ as well. Thus, $\Ext_R^i(M,R)=\Ext_R^i(\operatorname{Tr}(M),R)=0$ for all $i>0$. Therefore, the result follows from \cite[Proposition 2.8]{AuslanderReitenConjectureforNormalRings}.
\end{proof}
\begin{corollary}\label{432c}
    Let $R$ be a local ring of depth $t\geq 1$. Let $M$ be a locally free $R$-module on the punctured spectrum such that $\operatorname{G-dim}_R(M^\ast)<\infty$. If $\depth(M)\geq t$ and $\Ext_R^{t-1}(M,M)=0$, then $M$ is free. 
\end{corollary}
\begin{proof}
    This follows from Proposition \ref{r61} and Corollary \ref{a11p}. 
\end{proof}

\begin{theorem}\label{rpsCM}
    Let $R$ be a ring satisfying $(S_1)$ and $X$ be a subset of  $\Spec(R)$ containing $X^1(R)$. Let 
    \begin{center}
        $s:=\inf \left\{\operatorname{depth} R_{\mathfrak{p}} \mid \mathfrak{p} \in \operatorname{Spec}(R) \backslash X\right\}$ and $t:=\sup \left\{\operatorname{depth} R_{\mathfrak{p}} \mid \mathfrak{p} \in \operatorname{Spec}(R) \backslash X\right\}$. 
    \end{center}
    Let  $M$ be a locally free  $R$-module on $X$ such that $\operatorname{G-dim}_R(M^\ast)<\infty$. If $\Ext_R^i(M,R)=\Ext_R^j(M,M)=0$ for all $1\leq i \leq t$ and $s-1 \leq j \leq t-1$, then $M$ is projective.
\end{theorem}
\begin{proof}
    We may suppose that $X\not=\Spec(R)$. We show that $M_\mathfrak{p}$ is free for all $\mathfrak{p} \in \Spec(R)$. By assumption this is true if  $\mathfrak{p} \in X$. So, suppose that $\mathfrak{p} \in \Spec(R)\backslash X$. Therefore $\operatorname{ht}(\mathfrak{p})\geq 2$ and $\operatorname{depth}_{R_\mathfrak{p}}(R_\mathfrak{p})\geq 1$ since $R$ satisfies $(S_1)$. By definition of $s$ and $t$, note that $s\leq \operatorname{depth}_{R_\mathfrak{p}}(R_\mathfrak{p})\leq t$, so $\Ext_{R_\mathfrak{p}}^i(M_\mathfrak{p}, R_\mathfrak{p})=0$ for all $1\leq i \leq \operatorname{depth}_{R_\mathfrak{p}}(R_\mathfrak{p})$ and $\Ext_{R_\mathfrak{p}}^{\operatorname{depth}_{R_\mathfrak{p}}(R_\mathfrak{p})-1}(M_\mathfrak{p}, M_\mathfrak{p})=0$. Hence, by Proposition \ref{r61}, $M_\mathfrak{p}$ is free. 
\end{proof}

Observe that Proposition \ref{r61},  Corollary \ref{432c} and Theorem \ref{rpsCM} are generalizations of Theorem \ref{araq}.

Now, we explore some consequences of Theorem \ref{rpsCM} for Cohen-Macaulay local rings.
\begin{remark}(\cite[Remark 2.11]{AuslanderReitenConjectureforNormalRings}) \label{qas11t} Let $R$ be a Cohen-Macaulay local ring of dimension $d$ and $n$ be an integer such that $1\leq n \leq d-1$.  Then $$\inf \{\operatorname{depth}_{R_\mathfrak{p}}( R_{\mathfrak{p}}) \mid \mathfrak{p} \in \operatorname{Spec} (R) \backslash \mathrm{X}^n(R)\}=n+1\,\, \mbox{ and }\,\,$$ $$\sup \{\operatorname{depth}_{R_\mathfrak{p}}( R_{\mathfrak{p}}) \mid \mathfrak{p} \in \operatorname{Spec} (R) \backslash \mathrm{X}^n(R)\}=d.$$
\end{remark}
\begin{theorem}\label{exrt}
    Let $R$ be a Cohen-Macaulay local ring of dimension $d\geq 2$, $M$ be an $R$-module and $n$ be a positive integer. Suppose that $M$  is  locally of finite projective dimension on $X^n(R)$ and that $\Ext_R^j(M,M)=0$ for all $n\leq j \leq d-1$. Then the following hold:
    \begin{enumerate}
        \item If $\operatorname{G-dim}_R(M^\ast)<\infty$ and $\Ext_R^i(M,R)=0$ for all $1\leq i \leq d$, then $M$ is free.
        \item If $\operatorname{G-dim}_R(M)<n,$ then $\operatorname{pd}_R(M)<\infty$.
    \end{enumerate}
\end{theorem}
\begin{proof}
(1) We claim that $M$ is locally free on $X^{n}(R)$. In fact, given $\mathfrak{p} \in X^{n}(R)$, we have by hypothesis that $\operatorname{pd}_{R_\mathfrak{p}}(M_\mathfrak{p})<\infty$ and $\Ext_{R_\mathfrak{p}}^i(M_\mathfrak{p}, R_\mathfrak{p})=0$ for all $1\leq i \leq d$. Hence, by \cite[p. 154, Lemma 1(iii)]{CommutativeRingTheory}, $M_\mathfrak{p}$ is free as $R_\mathfrak{p}$-module. Thus item (1) follows immediately from Theorem \ref{rpsCM} and Remark \ref{qas11t}.

(2) Let $p=\operatorname{G-dim}_R(M)$. We prove the result by considering two cases.
 
  (2.a) Suppose $p=0$. In this case, $M$ is totally reflexive.  Hence $M^\ast$ is totally reflexive and $\Ext_R^i(M,R)=0$ for all $ i \geq 1$. So, it follows from (1) that $M$ is free. In particular, $\operatorname{pd}_R(M)<\infty$.

  (2.b) Suppose $p>0$. We have $\Ext_R^i(M,R)=0$ for all $i\geq p+1$. By \cite[Lemma 2.6]{AuslanderReitenConjectureforNormalRings}, we get $$\Ext_R^i(M,M)\cong \Ext_R^i\left(\Omega^p(M), \Omega^p (M)\right)$$ for all $i\geq n \geq p+1$. Hence $\Ext_R^i\left(\Omega^p (M), \Omega^p (M)\right)=0$ for all $n \leq i \leq d-1$. Note that $\operatorname{G-dim}_R(\Omega^p(M))=0$ and that $\Omega^p(M)$ is locally of finite projective dimension on $X^n(R)$. Hence, by (1), $\Omega^p(M)$ is free and therefore $\operatorname{pd}_R(M)<\infty$. 
\end{proof}

\begin{corollary}\label{ij}
    Let $R$ be a Cohen-Macaulay normal local ring of dimension $d\geq 2$, and let $M$ be an $R$-module such that $\operatorname{G-dim}_R(M^\ast)<\infty$. Suppose that $\Ext_R^i(M,R)=\Ext_R^j(M,M)=0$ for all $1\leq i \leq d$ and  $1\leq j \leq d-1$. Then $M$ is free. 
\end{corollary}

\subsection*{On the Auslander-Reiten conjecture, and a result of Dey and Ghosh} 
In \cite[Corollary 6.9(2)]{FiniteHomologicalDimensionOfHomv3}, Dey and Ghosh demonstrated that the Auslander-Reiten conjecture holds for (finitely generated)  $R$-modules $M$ satisfying $\operatorname{G-dim}_R(M)<\infty$ and \linebreak $\operatorname{pd}_R(\Hom_R(M,M))<\infty$. Motivated by this, we investigate whether the conjecture remains valid when we consider $\operatorname{G-dim}_R(M^\ast)<\infty$ instead of $\operatorname{G-dim}_R(M)<\infty$ in the referenced case.

\begin{theorem}\label{fjk}
     Let $R$ be a local ring of depth $t$, and let $M$  be an $R$-module such that $\operatorname{G-dim}_R(M^\ast)<\infty$. Then the following conditions are equivalent. 
     \begin{enumerate}[(1)]
         \item $M$ is free.
         \item $\Hom_R(M,M)$ is free and $\Ext_R^j(M,M)=0$ for all $1\leq j \leq t-1$.
         \item  $\Hom_R(M,M)$ has finite projective dimension, and $\Ext_R^i(M,R)= \Ext_R^j(M,M) =0$ for all $1\leq i \leq t$ and $1\leq j \leq t-1$. 
     \end{enumerate}
\end{theorem}
\begin{proof}
    $(1) \Rightarrow (2)$ is trivial. 

    $(2) \Rightarrow (3)$. We only need to prove that $\Ext_R^i(M,R)=0$ for all $1\leq i \leq t$. For this, we will show that $M$ is totally reflexive. In view of \cite[Corollary 7.6(2)]{FiniteHomologicalDimensionOfHomv3}, $M$ satisfies $(\widetilde{S}_t)$. Therefore, $M$ is totally reflexive by Theorem \ref{pos5}.

    $(3) \Rightarrow (1)$. By Theorem \ref{pos5}, $M$ is totally reflexive. Therefore by \cite[Corollary 6.8]{FiniteHomologicalDimensionOfHomv3}, $M$ is free. 
\end{proof}

\begin{corollary} Let $R$ be a ring. The Auslander-Reiten conjecture holds true for all (finitely generated) $R$-modules $M$ such that $\operatorname{G-dim}_R(M^\ast)<\infty$ and $\operatorname{pd}_R( \Hom_R(M,M))<\infty$.
\end{corollary}

An $R$-module $C$ is said to be \textit{semidualizing} if the natural map $R \rightarrow \Hom_R(C,C)$ is an isomorphism and $\Ext_R^i(C,C)=0$ for all $i>0$. 

\begin{corollary}\label{sd}
    Let $R$ be a local ring and let $C$ be a semidualizing $R$-module. If 
 $\operatorname{G-dim}_R(C^\ast)<\infty$, then $C\cong R$.
    \end{corollary}

\section{Gorenstein criteria and related questions}
Throughout this section, let $(R, \mathfrak{m},k)$ be a local ring. In this section, we aim to provide a number of criteria for a local ring $R$ to be Gorenstein in terms of the dual of certain modules having finite Gorenstein dimension. In this sense, we discuss some related questions to this subject.

The first criterion that we present in this section is an application of  Corollary \ref{sdq}. Before stating it, first recall that if $R$ is a ring with total quotient ring $Q$, then an  $R$-module $M$ is said to have a (generic) \textit{rank}, denoted by $\operatorname{rank}(M)$, if $M \otimes_R Q$ is a free $Q$-module of rank $\operatorname{rank}(M)$. Let $e(R)$ denote the Hilbert-Samuel multiplicity of $R$ and $\mu(M)$ denote the minimum number of generators of $M$. Moreover, if $M$ is Cohen-Macaulay and $e(M) = \mu(M)$, then $M$ it is said to be an \textit{Ulrich} module (\cite[Definition 2.1]{AlmostGorensteinRingsTowardsATheoryOfHigherDimension}).

 \begin{proposition}\label{Gorentein2}  
Let $R$ be a Cohen-Macaulay local ring, and let 
$M$ be an $R$-module such that $\operatorname{G-dim}_R(M^\ast)<\infty$.  Then $R$ is Gorenstein in each one of the following cases:
\begin{enumerate}[(1)]
    \item $M$ is Cohen-Macaulay with positive rank and $2\mu(M) > e(R) \operatorname{rank}(M)$. 
    \item $M$ is an Ulrich module and $M^\ast \not=0$. 
\end{enumerate}
 \end{proposition}
 \begin{proof}
(1) The positivity of the rank of $M$ implies that $M^\ast \neq 0$. 

Since $M$ is Cohen-Macaulay and $\operatorname{G-dim}_R(M^\ast) < \infty$, according to Corollary \ref{sdq}, we conclude that $M$ is totally reflexive. Consequently, $\Ext_R^i(M,R) = 0$ for all $1 \leq i \leq d$. Thus, by assumption and \cite[Theorem 3.1]{Ulrich}, we obtain the result.

(2) By definition, $M$ is Cohen-Macaulay, and by Corollary \ref{sdq}, $M$ is totally reflexive. Hence $R$ is Gorenstein by \cite[Proposition 2.19]{ComplexityandrigidityofUlrichModules}. 
\end{proof}

Due to Proposition \ref{Gorentein2}, it is natural to ask the following. 

\begin{question} \label{pow} Let $R$ be a local ring. Suppose there exists a non-free (finitely generated) $R$-module $M$ such that $M^\ast\not=0$ and $\operatorname{G-dim}_R(M^\ast) < \infty$. Then, is $R$ Gorenstein?
\end{question}

\begin{remark} This question is a refined version of Question 4.5 of the second version of this paper,  which was answered negatively in \cite[Example 3.10 and Remark 3.11]{ARCconjectureForModulesWhoseSelfDualHasFiniteCompleteIntersectionDimension} by Ghosh and Samanta. The same example shows that Question \ref{pow} is false in general. 
\end{remark}

On the other hand, it is known that if $R$ is a local ring such that $\operatorname{G-dim}_R(M)<\infty$ for all $R$-modules, then $R$ is Gorenstein (\cite[Theorem 1.4.9]{GorensteinDimensions}). Thus, by considering the dual of $R$-modules it is natural to ask the following. 

\begin{question}\label{GorA} Let $R$ be a local ring. 
If $\operatorname{G-dim}_R(M^\ast)<\infty$ for all (finitely generated) $R$-modules $M$, then is $R$ Gorenstein?
\end{question}
Then motivated by Question \ref{GorA}, we obtain the following result which says that its answer is positive. 
\begin{theorem}\label{GG} Let $R$ be a local ring. If $\operatorname{G-dim}_R(M^\ast)<\infty$ for all (finitely generated) $R$-modules $M$, then $R$ is Gorenstein. 
\end{theorem}
\begin{proof} First, suppose that $\operatorname{depth}_R(R)=0$. Then $k^\ast \not=0$ and by assumption \break
$\operatorname{G-dim}_R(k^\ast)<\infty$. Since $k^\ast \cong k^n$ for some $n\geq 1$, note that $\operatorname{G-dim}_R(k)<\infty$. Hence $R$ is Gorenstein.

Now, suppose that  $\operatorname{depth}_R(R)\geq 1$. We may write $\mathfrak{m}=(x_1, \ldots, x_n)$. Let $\varphi: R \to R^n$ be the homomorphism defined by $\varphi(r)=(rx_1, \ldots, rx_n)$. Since $\operatorname{depth}_R(R)>0$, we see that $\varphi$ is injective. Thus we have an exact sequence 
    $$0 \rightarrow R \stackrel{\varphi}{\rightarrow} R^{n} \rightarrow M \rightarrow 0.$$
  Dualizing this exact sequence we obtain a sequence 
    \begin{equation}\label{eqsi}
        0 \rightarrow M^\ast \to R^n \stackrel{\varphi^\ast }{\rightarrow} R \to \operatorname{coker}(\varphi^\ast)=k \to 0.
    \end{equation} 
Since $\operatorname{G-dim}_R(M^\ast)<\infty$, from \eqref{eqsi}, we see that $\operatorname{G-dim}_R(k)<\infty,$ concluding that $R$ is Gorenstein. 
\end{proof}
\begin{remark} Theorem \ref{GG} can be obtained in other ways. For instance, it could be derived from \cite[Corollary 3.5]{ktorsionlessmoduleswithfiniteGorensteindimension}. Additionally, the theorem can be obtained as an application of Theorem \ref{GteosPa}. Indeed, setting $d=\dim_R(R)$, by the assumption and Theorem \ref{GteosPa}, all (finitely generated) $R$-modules satisfy the following equivalences:
    \begin{center}
      $M$ is $(d+1)$-torsionfree $\Longleftrightarrow$ $M$ is $(d+1)$-syzygy $\Longleftrightarrow$ $M$ satisfies $(\widetilde{S}_{d+1}).$
    \end{center}
   Then by \cite[Theorem 1.4]{WhenArenSysymodulesnTorsionfree}, $R$ satisfies $(G_{d})$, concluding that $R$ is Gorenstein.    
\end{remark}

For a local ring of depth at most one, the finiteness of the Gorenstein dimension of the dual of its maximal ideal characterizes its Gorensteiness.

\begin{proposition}\label{sq884}
    Let $R$ be a local ring of depth at most one. If $\operatorname{G-dim}_R(\mathfrak{m}^\ast)<\infty$, then $R$ is Gorenstein. 
\end{proposition}
\begin{proof}
    First, suppose $\operatorname{depth}_R(R)=0$. Since $\operatorname{G-dim}_R(\mathfrak{m}^\ast)<\infty, $ then $\operatorname{G-dim}_R(\operatorname{Tr} (\mathfrak{m}))<\infty$. Since $\operatorname{depth}_R(R)=0$, the Auslander-Bridger formula shows that $\operatorname{Tr} (\mathfrak{m})$ is totally reflexive and hence $\mathfrak{m}$ as well. Now, by considering the exact sequence $0 \to \mathfrak{m} \to R \to k \to 0,$ we see that $\operatorname{G-dim}_R(k)<\infty$. Thus, we conclude that $R$ is Gorenstein.

    Now, suppose that $\operatorname{depth}_R(R)=1$. Then the exact sequence $0 \to \mathfrak{m} \to R \to k \to 0,$ induces an exact sequence 
    $$0 \to k^\ast \to R \to \mathfrak{m}^\ast \to \Ext_R^1(k,R) \to 0.$$
    Since $\operatorname{depth}_R(R)=1$, we see that $k^\ast=0$ and $\Ext_R^1(k,R)\cong k^n$ for some $n\geq 1$. Thus, we obtain an exact sequence 
    $$0 \to R \to \mathfrak{m}^\ast \to k^n \to 0.$$
    In view of that $\operatorname{G-dim}_R(\mathfrak{m}^\ast)<\infty, $ we see that $\operatorname{G-dim}_R(k)<\infty$. This shows that $R$ is Gorenstein. 
\end{proof}

Motivated by Question \ref{sqeeww}, we consider the following condition:

\begin{description}
    \item[(GDUAL)] Every (finitely generated) $R$-module $M$ whose dual $M^\ast$ is of finite Gorenstein dimension, is also of finite Gorenstein dimension.
\end{description}

\begin{remark} If the local ring $R$ satisfies $\operatorname{(GDUAL)}$, then $R$ is not necessarily Gorenstein unless $\operatorname{depth}_R(R)>0$. Indeed, suppose $\operatorname{depth}_R(R)=0$. From the Auslander-Bridger formula, we can see that all $R$-modules of finite Gorenstein dimension are totally reflexive. Thus, any $R$-module whose dual has finite Gorenstein dimension, is totally reflexive (see Remark \ref{wss} and Facts \ref{eqpl}(2)), and hence $R$ satisfies $\operatorname{(GDUAL)}$. Hence the existence of non-Gorenstein local rings of depth zero shows that $R$ is not necessarily Gorenstein. The case $\depth_R(R)>0$ is discussed in the next proposition. 
\end{remark}

Motivated by condition $(\operatorname{GDUAL})$, we define a weaker condition than it, and prove that if $R$ is of positive depth and satisfies such condition, then $R$ is Gorenstein. 
\begin{description}
    \item[(GDUAL*)] Every (finitely generated) $R$-module $M$ whose dual $M^\ast$ is non-zero and of finite Gorenstein dimension, is also of finite Gorenstein dimension.
\end{description}
 
\begin{proposition}\label{Gdual}
Let $R$ be a local ring of  depth at least one satisfying $\operatorname{(GDUAL*)}$. Then $R$ is Gorenstein.
\end{proposition}
\begin{proof}
Since $R$ has positive depth, note that  $k^\ast=0$. Let $M=R \oplus k$. Then $M^\ast=R$ and hence $\operatorname{G-dim}_R(M^\ast)<\infty$. Thus $\operatorname{G-dim}_R(k)<\infty$ by assumption. Hence, $R$ is Gorenstein. 
\end{proof}

Now, motivated by  Corollary \ref{sd}, we will discuss about a question that was addressed by Holanda and Miranda-Neto in \cite{VanishingOf(Co)homologyfrenessCriteriaAndTheAuslanderReitenConjectureforCohen-MacaulayBurch}. The question is as follows:

\begin{question}\label{qp5}(\cite[Question 5.24]{VanishingOf(Co)homologyfrenessCriteriaAndTheAuslanderReitenConjectureforCohen-MacaulayBurch}) Let $R$ be a Cohen-Macaulay local ring with a canonical module $\omega$. If $\operatorname{pd}_R( \omega^\ast)<\infty$, then must $R$ be Gorenstein?    
\end{question}

Firstly, we observe that the answer to Question \ref{qp5} is positive and follows directly from the Corollary \ref{sd}. However, it is worth highlighting that this question had already been answered affirmatively, albeit through a different proof, by Asgharzadeh in \cite[Theorem 11.3]{ReflexivityRevisited}. Additionally, it is worth noting that long before Holanda and Miranda-Neto posed this question, the same result could be obtained directly from the Foxby equivalence (see \cite[Theorem 3.4.11]{GorensteinDimensions}).

We finish this section by proposing the following question, which is motivated by a celebrated result of Foxby \cite{FoxbyIsomorphimsbetween} as well as a result of Holm \cite{RingswithfiniteGorensteininjectivedimension}. The first (resp. second) asserts that the existence of a non-zero (finitely generated) $R$-module $M$ of projective dimension (resp. Gorenstein dimension) and injective dimension both finite implies that $R$ is Gorenstein. 
\begin{question}\label{qw1}
If there exists a (finitely generated) $R$-module $M$ such that\break $\operatorname{G-dim}_R(M^\ast)<\infty$ and $\operatorname{id}_R(M)<\infty$, then is $R$ Gorenstein?   
\end{question}

\section{Preliminary facts on K\"ahler differentials}\label{section7}
In this section, we introduce some notations and recall known results on Kähler and differential modules in the context of algebras and locally ringed spaces for convenience.
\subsection*{Affine derivation modules and K\"ahler differentials $n$-th order} 
Let  $S$ be a (commutative) ring and  $R$ be an $S$-algebra and $M$ be an $R$-module. Recall that an $S$-linear map $D : R \rightarrow M$ is said to be a \textit{derivation} if for any two elements $x_0, x_1 \in R$, the following identity holds:
\[ D(x_0 x_1) = x_0 D(x_1) + x_1 D(x_0). \]
A derivation of order $n$ can be defined generalizing the previous identity as follows.

 An $S$-linear map $D : R \rightarrow M$ is said to be an $n$-{\it th order derivation} if for any $x_0, \ldots, x_n \in R$, the following identity holds:
\[ D(x_0 \cdots x_n) = \sum_{s=1}^{n} (-1)^{s-1} \sum_{i_1 < \cdots < i_s} x_{i_1} \cdots x_{i_s} D(x_0 \ldots \hat{x}_{i_1} \cdots \hat{x}_{i_s} \cdots x_n), \]
where $\hat{x}_{i_j}$ means that this element does not appear in the product.  The set of $n$-th order derivations of an $S$-algebra $R$ into an $R$-module $M$ over $S$ will be denoted by $\text{Der}^n_S(R, M)$. When $M = R$, we shall use the notation $\text{Der}^n_S(R)$ in place of $\text{Der}^n_S(R, R)$.

The module of derivations of order $n$ can be represented as follows. Let $I$ denote the kernel of the homomorphism $R \otimes_S R \rightarrow R$, $a \otimes b \mapsto ab$. Giving structure of $R$-module to $R \otimes_S R$ by multiplying on the left, we define the $R$-module
$$\Omega^{(n)}_{R/S} := I/I^{n+1}.$$

Define the map $d_n^R : R \rightarrow \Omega^{(n)}_{R/S}$,\, $a \mapsto (1 \otimes a - a \otimes 1) + I^{n+1}$. This map is a derivation of order $n$, and its image generates $\Omega^{(n)}_{R/S}$ as an $R$-module (see \cite[Chapter II-1]{Nakai}).

 The $R$-module $\Omega^{(n)}_{R/S}$ is called the {\it module of K\"ahler differentials of order} $n$ of $R$ over $S$. The map $d_n^R$ is called the \textit{canonical derivation} of $R$ in $\Omega^{(n)}_{R/S}$. It comes equipped with a universal derivation $d_{S/R}\in \text{Der}_R^n\left(S, \Omega^{(n)}_{S/R}\right)$ with the property that composition with $d_{S/R}$ yields an isomorphism $\Hom_R\left(\Omega^{(n)}_{R/S},R\right) \cong \text{Der}_S^n(R)$ (\cite[Proposition 1.6]{O}). To see more properties regarding modules of derivations and K\"ahler differentials, we recommend \cite{O, NKI, Nakai}.

\begin{remark}\label{finites}
    It is important to note that the modules of derivations and the K\"ahler differentials may not be finitely generated. However, they are finitely generated in certain cases. For instance: If $R$ is essentially of finite type over $S$.
If $S=k$ is a field with a valuation and $R$ is an analytic $k$-algebra, meaning $R$ is module-finite over a convergent power series ring $k\{x_1,\ldots, x_n\}$.
 If $S=k$ is a field, $(R, \mathfrak{m})$ is a complete local ring, and $R/\mathfrak{m}$ is a finite extension of $k$. For this reason, in this article, for each $n$, we will consider $\Omega^{(n)}_{R/S}$ to be finitely generated as an $R$-module. In particular, by the universal property, $\text{Der}^n_S(R)$ is also finitely generated as an $R$-module.
\end{remark}

We now present key definitions and properties that are central to this section.

\begin{definition}
 A {\it scheme} is a ringed topological space $(X,\mathcal{O}_X)$ admitting an open covering $\{U_i\}_i$ such that $(U_i, \mathcal{O}_X|_{U_i} )$ is an affine scheme for every $i$. Therefore, a scheme is a locally ringed space.  A scheme $(X,\mathcal{O}_X)$ is said to be {\it Noetherian} if it is a finite union of affine open $\Spec (R_i)$ such that $R_i$ is a (commutative) Noetherian ring for every $i$. We say that a scheme is {\it locally Noetherian} if every point has a Noetherian open neighborhood. 
\end{definition}
We can see that if  $(X,\mathcal{O}_X)$ is a locally Noetherian scheme, then the local rings \( \mathcal{O}_{X,x} \) are Noetherian rings. Often we simply write $X$ instead of  $(X,\mathcal{O}_X)$.

\subsection{Global K\"ahler differentials}
Next, we recall the global definition of the constructions given above. To do this, instead of gluing together the modules $\Omega^{(1)}_{S/R}$, we provide a global definition, and then point out that it reduces to the original definition given above. Before defining it, it is worth noting that the fiber product exists in our context (see \cite[Theorem 4.18]{GW}).
Let $\varphi:X\to Y$ be a morphism of schemes.  This induces a morphism $\Delta:X\to X\times_Y X$, called the  diagonal morphism.
 The map $\Delta$ comes with a map of sheaves $\Delta^{\#}: \mathcal{O}_{X\times_YX}\to \Delta_{*} \mathcal{O}_X$ on $X\times_YX$. Although really a map of sheaves of rings, we will regard it as a map of sheaves of $\mathcal{O}_{X\times_YX}$-modules. Let $\mathcal{I}$ be its kernel, again regarded as an $\mathcal{O}_{X\times_YX}$-module.

\begin{definition} The {\it sheaf relative K\"ahler differential} of $\mathcal{O}_X$-modules is defined as:
\[
\Omega^{(1)}_{X/Y} = \Delta^{\ast}(\mathcal{I}/\mathcal{I}^2).
\]
\end{definition}

 Thus, more generally, the {\it sheaf of relative K\"ahler $n$-differentials} of $X$ over $Y$ is defined by
$$
\Omega_{X/Y}^{(n)}:=
\begin{cases}
      \bigwedge^n \Omega_{X/Y}^{(1)}, & \text{if } n\geq 1,\\
      \mathcal{O}_X, & \text{if } n=0,\\
      0, & \text{if } n< 0.
\end{cases}
$$
Note that if $Y=\Spec(k)$, where $k$ is a field, $\Omega_{X/Y}^{(n)}$ is usually denoted by $\Omega_{X}^{(n)}$. For further details and properties, see for example, \cite[Section I.6.6, p. 34]{Gro1}.

\begin{facts}\label{fact1}
Let $\varphi: X \to Y$ be a morphism of locally Noetherian schemes. 
\begin{enumerate}
\item \mbox{(\cite[Proposition 6.1.20]{Qing} or \cite[Corollaire (I6.4.22)]{Gro1})}. If $\varphi:X\to Y$ is locally of finite type, then $\Omega_{X/Y}^{(1)}$ is coherent. Additionally, by \cite[Proposition 7.48]{GW}, $\Omega_{X/Y}^{(n)}$ is also a coherent $\mathcal{O}_X$-module  for all $n>1$. 

\item \mbox{ (\cite[Proposition 1.17, ch. 6]{Qing}).} For each $x \in X$, there is a canonical isomorphism of $\mathcal{O}_{X,x}$-modules
$$\left(\Omega_{X/Y}^{(1)}\right)_x=\Omega_{\mathcal{O}_{X,x}/\mathcal{O}_{Y,\varphi(x)}}^{(1)}.$$ Since the stalk commutes with exterior power, more generally, one obtains $\left(\Omega_{X/Y}^{(n)}\right)_x\cong \Omega_{\mathcal{O}_{X,x}/\mathcal{O}_{Y,\varphi(x)}}^{(n)}=: \Omega_{(X,x)/(Y,\varphi(x))}^{(n)}$.

\item \mbox{(\cite[Corollaire (I6.5.5)]{Gro1})}. There is also a definition of the sheaf of $Y$-derivations over $X$, denoted by $\mathcal{D}er_Y^{1}(X)$. Moreover, the universal property of $\Omega_{X/Y}^{(1)}$ comes down to saying that there is an isomorphism of $\mathcal{O}_X$-modules
$\mathcal{D}er_Y^{1}(X)\cong \mathcal{H}om_{\mathcal{O}_X}\left(\Omega_{X/Y}^{(1)},\mathcal{O}_X\right).$ Naturally, is defined the sheaf of $n$-th order $Y$-derivations over $X$.
 $$\mathcal{D}er_Y^{n}(X):=\mathcal{H}om_{\mathcal{O}_X}\left(\Omega_{X/Y}^{(n)},\mathcal{O}_X\right).$$

  \item   \mbox{(\cite[Lemma 28.20.2, \href{https://stacks.math.columbia.edu/tag/05P1}{Tag 05P1}]{stacks-project}).} Let $X$ be a locally Noetherian scheme, and let $\mathcal{F}$ be a coherent $\mathcal{O}_X$-module. Then  $\mathcal{F}$ is locally free if and only if for all $x \in X$, the stalk $\mathcal{F}_x$ is a free $\mathcal{O}_{X,x}$-module.
\end{enumerate}
\end{facts}

We provide an elementary lemma that may be familiar, but we have not found it in the literature.
\begin{lemma}\label{lemaseup} 
     Let \(X\) and \(Y\) be locally Noetherian schemes, and let \(\varphi: X \to Y\) be a morphism locally of finite type. Then, the following properties hold:
\begin{enumerate}
    \item[(i)] $\mathcal{D}er_Y^{n}(X)$ is a coherent $\mathcal{O}_X$-module.
    \item[(ii)] For $x \in X$ there is a canonical isomorphism of $\mathcal{O}_{X,x}$-modules $$(\mathcal{D}er_Y^{n}(X))_x\cong {\rm Der}^n_{\mathcal{O}_{Y,\varphi(x)}}(\mathcal{O}_{X,x}).$$ 
    \end{enumerate}
\end{lemma}
\begin{proof} 
(i) From Facts \ref{fact1}(1), we establish the coherence of $\Omega_{X/Y}^{(n)}$. Subsequently, from  the equality $\mathcal{D}er_Y^{n}(X)=\mathcal{H}om_{\mathcal{O}_X}\left(\Omega_{X/Y}^{(n)},\mathcal{O}_X\right)$ and \cite[Exercise 5.1.6(b)]{Qing}, we deduce the result.

(ii) 
\begin{align*}
     (\mathcal{D}er_Y^{n}(X))_x & \cong \left(\mathcal{H}om_{\mathcal{O}_X}\left(\Omega_{X/Y}^{(n)},\mathcal{O}_X\right)\right)_x\,\,\, \text{from Facts \ref{fact1}(3)}  \\
     & \cong \Hom_{\mathcal{O}_{X,x}}\left(\Omega_{\mathcal{O}_{X,x}/\mathcal{O}_{Y,\varphi(x)}}^{(n)},\mathcal{O}_{X,x}\right) \\
     &\cong {\rm Der}^n_{\mathcal{O}_{Y,\varphi(x)}}(\mathcal{O}_{X,x}).
\end{align*}
Since  $\Omega_{X/Y}^{(n)}$ is coherent $\mathcal{O}_X$-module (Facts \ref{fact1}(1)), so $\mathcal{H}om$ commutes with taking stalks 
(\cite[Proposition 7.27]{GW}). Now from Facts \ref{fact1}(2)(3) from  the universal property ( \cite[Theorem 2.2.6]{HS} or \cite[Proposition 1.6]{O}) we get the last isomorphism.
\end{proof}

Based on Remarks \ref{finites}, Facts \ref{fact1}, and Lemma \ref{lemaseup}, the schemes considered will be locally Noetherian, and the sheaves \(\mathcal{D}er_Y^{n}(X)\) and \(\Omega_{X/Y}^{(n)}\) will always be non-zero coherent sheaves. Furthermore, we assume that these sheaves have non-zero stalks, implying that the chosen elements lie in the support.
For simplicity, the stalks of \(\Omega_{(X,x)/(Y,\varphi(x))}^{(n)}\) and \({\rm Der}^n_{\mathcal{O}_{Y,\varphi(x)}}(\mathcal{O}_{X,x})\) will be denoted by \(\Omega_{(X,x)/(Y,y)}^{(n)}\) and \({\rm Der}^n_{(Y,y)}(X,x)\), respectively, where \(y = \varphi(x)\).

\section{Applications}\label{section8}

\subsection{Applications on the $k$-torsion of the modules of differentials}

Next, we study the conditions under which the Kähler differential module $\Omega_{X/Y}^{(n)}$ over schemes is reflexive, $k$-torsionfree, and $k$-syzygy, with particular attention to the cases where $k = 1$ and $k = 2$. The investigation of these properties in Kähler differential modules of order $n$ concerning regularity questions of affine rings has a long history in the literature, addressed by various authors, especially in the context of algebraic varieties and analytic spaces, particularly when $Y$ is the spectrum of field, see, for example, \cite{Kunz, Milher, graf, Lipman, Suzuki, Vetter, Lebelt, Greuel}). Recently, Graf and Milher in \cite{graf, Milher} gave new results on the connection between torsion, cotorsion, and reflexivity of the Kähler differential modules of order $n$ in complex analytic varieties and algebraic varieties. For a summary of the existing results, see \cite[p. 133]{Milher}. Motivated by these results, we provide the following results in the context of schemes.

\begin{theorem}\label{stalksquema1}
    Let \(X\) and \(Y\) be locally Noetherian schemes, $x\in X$ and $k$ be a non-negative integer such that ${\rm Der}^n_{{(Y,y)}}(X,x)$ has locally finite Gorenstein dimension on $\widetilde{X}^{k-1}(\mathcal{O}_{X,x})$. Then the following conditions are equivalent:
    \begin{enumerate}
        \item $\Omega_{(X,x)/(Y,y)}^{(n)}$ is $k$-torsionfree.
        \item $\Omega_{(X,x)/(Y,y)}^{(n)}$ is $k$-syzygy.
        \item $\Omega_{(X,x)/(Y,y)}^{(n)}$ satisfies $(\widetilde{S}_k)$. 
    \end{enumerate}
\end{theorem}
\begin{proof}
    By Facts \ref{fact1}(2), we have $\Omega_{(X,x)/(Y,y)}^{(n)}=\Omega_{\mathcal{O}_{X,x}/\mathcal{O}_{Y,y}}^{(n)}$. Therefore, $\Omega_{(X,x)/(Y,y)}^{(n)}$ is an $\mathcal{O}_{X,x}$-module. Moreover, we know that ${\rm Der}^n_{{(Y,y)}}(X,x)=\left(\Omega_{(X,x)/(Y,y)}^{(n)}\right)^{\ast}$. Now, the required equivalences follow from Theorem \ref{GteosPa}.
\end{proof}

As an immediate consequence of Theorem \ref{stalksquema1}, or Corollaries \ref{saq54} and \ref{saq55}, we have the following corollaries.

\begin{corollary}\label{corsche1} 
    Let \(X\) and \(Y\) be locally Noetherian schemes, and let $x\in X$ such that ${\rm Der}^n_{{(Y,y)}}(X,x)$ has locally finite Gorenstein dimension on ${\rm Ass}_{\mathcal{O}_{X,x}}(\mathcal{O}_{X,x})$. Then the following conditions are equivalent:
    \begin{enumerate}
        \item $\Omega_{(X,x)/(Y,y)}^{(n)}$ is $1$-torsionfree.
        \item $\Omega_{(X,x)/(Y,y)}^{(n)}$ is $1$-syzygy.
        \item $\Omega_{(X,x)/(Y,y)}^{(n)}$ satisfies $(\widetilde{S}_1)$. 
    \end{enumerate}
\end{corollary}

\begin{corollary}\label{corsche2}    
Let \(X\) and \(Y\) be locally Noetherian schemes, and let \(x\in X\) such that \({\rm Der}^n_{{(Y,y)}}(X,x)\) has locally finite Gorenstein dimension on \(\widetilde{X}^{1}(\mathcal{O}_{X,x})\). Then the following conditions are equivalent:
    \begin{enumerate}
        \item $\Omega_{(X,x)/(Y,y)}^{(n)}$ is $2$-torsionfree.
        \item $\Omega_{(X,x)/(Y,y)}^{(n)}$ is $2$-syzygy.
        \item $\Omega_{(X,x)/(Y,y)}^{(n)}$ satisfies $(\widetilde{S}_2)$.
        \item  $\Omega_{(X,x)/(Y,y)}^{(n)}$ is reflexive.
    \end{enumerate}
\end{corollary}

\begin{remark}\label{afim}
    Let $\varphi : X \rightarrow Y$ be a morphism of schemes which is locally of finite type. If we consider $X$ and $Y$ as affine Noetherian schemes, that is, $X=\Spec(R)$ and $Y=\Spec(S)$ where $R$ and $S$ are local Noetherian rings. By Facts \ref{fact1}(2) and Lemma \ref{lemaseup}(ii), if we consider $x\in X$ and $\varphi(x)\in Y$ closed points respectively, then we have $\Omega_{(X,x)/(Y,\varphi(x))}^{(n)}=\Omega^{(n)}_{R/S}$ and ${\rm Der}^n_{{(Y,\varphi(x))}}(X,x)={\rm Der}^n_S(R)$. Moreover, these are finitely generated $R$-modules.
\end{remark}

From Remark \ref{afim}, Theorem \ref{stalksquema1}, and Corollaries \ref{corsche1} and \ref{corsche2}, we can obtain similar results in the affine context, that is, the derivation modules and Kähler differential modules defined over algebras.

\subsection{Applications on Herzog-Vasconcelos's conjecture}\label{section9}

One of the main conjectures, which may have motivated the study of properties of the K\"ahler differential modules and the regularity of schemes, is the famous Lipman-Zariski conjecture \cite{Lipman}, which states the following: Let \( X \) be a complex variety such that the tangent sheaf \( \mathcal{T}_X :=\mathcal{H}om_{\mathcal{O}_X}\left(\Omega_{X}^{(1)},\mathcal{O}_X\right) \) is locally free. Then \( X \) is smooth. This conjecture remains widely open. Lipman-Zariski conjecture, have
been extensively studied and resolved in special cases. References of Lipman-Zariski
conjecture can be found in Miller-Vassiliadou [37, Section 4]. It is important to note that this conjecture is generally false for positive characteristic. For example, consider the surface \( X \subset \mathbb{A}^3_k \) over a perfect field \( k \) of characteristic \( p \), defined by the equation \( xy - z^n = 0 \), where \( p \) divides \( n \). Then one can see that \( \mathcal{T}_X \) is free (see \cite[p. 892]{Lipman}).  Following the same direction as the Zariski-Lipman conjecture, there is also a homological version proposed by Herzog and Vasconcelos:

\begin{conjecture}[Herzog-Vasconcelos's conjecture]
Let $R$ be a local $k$-algebra, where $k$ is a field of characteristic zero. If ${\rm pd}_R({\rm Der}_k(R)) < \infty$, then ${\rm Der}_k(R)$ is a free $R$-module.
\end{conjecture}

 This conjecture, unlike the former, appears to be widely open (with the exception of a few specific cases; see, e.g., \cite[Section 4]{He-M} and \cite{HV, PerezNeto}).

Motivated by Herzog-Vasconcelos’s conjecture and the definitions from the previous sections, we can present the following question. But first, let us define: Let \( M \) be an \( R \)-module. Since
$
\operatorname{pd}_R(M) = \sup\{\operatorname{pd}_{R_{\mathfrak{p}}}(M_{\mathfrak{p}}) \mid \mathfrak{p} \in \Spec(R)\},
$
the {\it projective dimension of a coherent sheaf} $\mathcal{F}$ over an scheme \( X \) is defined as
$
\pd_{\mathcal{O}_X}\mathcal{F} = \sup\{\pd_{\mathcal{O}_{X,x}}(\mathcal{F}_x) \mid x \in X\}.
$

\begin{questions}\label{questionHV}({\bf Generalization of  Herzog-Vasconcelos's conjecture (GHVC)}).
\begin{itemize}
\item[(i)] Let $n\geq 1$ and  ${\rm pd}_{\mathcal{O}_X}{\mathcal{D}er}_Y^n(X) \, < \, \infty$. Under what assumptions on $X$ and $Y$, and for which values of $n$, does this imply that ${\mathcal{D}er}_Y^n(X)$ is locally free?
\item[(ii)] In particular, consider \(X = \operatorname{Spec}(R)\) and \(Y = \operatorname{Spec}(S)\) as in Remark \ref{afim}. Assume that for some integer $n\geq 1$, \( {\rm pd}_{R} (\operatorname{Der}_S^n(R)) < \infty \). Under what assumptions on \( R \), and for which values of \( n \) does this imply that \( \operatorname{Der}_S^n(R) \) is free?
\end{itemize}
\end{questions}

In the context of Question \ref{questionHV}(ii), a generalized form of the Zariski-Lipman Conjecture can be stated as follows: If \( k \) is a field and \( R \) is a \( k \)-algebra, then \( R \) is a regular local ring if \( \operatorname{Der}_k^n(R) \) is a free \( R \)-module for some \( n \geq 1 \). It should also be noted that for \( n > 1 \), Question \ref{questionHV}(ii) generally fails when \( S = k \) is a field of positive characteristic. For a counterexample, see \cite{Lu}.

A coherent sheaf $\mathcal{F}$ over a locally Noetherian scheme $X$ is said to be  \textit{Cohen-Macaulay} if, for every point $x \in X$, the stalk $\mathcal{F}_x$ is a Cohen-Macaulay module over the local ring $\mathcal{O}_{X,x}$.   

The following result provides a partial answer to the GHVC conjecture when \(X\) and $Y$ are locally Noetherian schemes with $X$ and the Kähler differential module \(\Omega_{X/Y}^{(n)}\) being Cohen-Macaulay.

\begin{proposition}\label{diff1}
Let \(X\) and \(Y\) be locally Noetherian schemes, with \(X\) Cohen-Macaulay, and suppose that \({\rm pd}_{\mathcal{O}_X}{\mathcal{D}er}_Y^n(X) < \infty\) for some $n\geq 1$. If $\Omega_{X/Y}^{(n)}$ is  Cohen-Macaulay, then \({\mathcal{D}er}_Y^n(X)\) is locally free.
 \end{proposition}
 \begin{proof} 
It is suffices to show it at the level of stalks. Note that, \({\rm pd}_{\mathcal{O}_X}{\mathcal{D}er}_Y^n(X) < \infty\) if and only if \(\operatorname{pd}_{\mathcal{O}_{X,x}}\left({\rm Der}^n_{\mathcal{O}_{Y,y}}(\mathcal{O}_{X,x})\right)<\infty\) for all \(x\in X\). Moreover, $\Omega_{X/Y}^{(n)}$ is  Cohen-Macaulay if and only if \(\Omega_{(X,x)/(Y,y)}^{(n)}\) is  Cohen-Macaulay for all \(x\in X\). Now, by Corollary \ref{erwq},  \(\Omega_{(X,x)/(Y,y)}^{(n)}\) is a free \(\mathcal{O}_{X,x}\)-module for all \(x\in X\). Thus, by the universal property, we obtain that \({\rm Der}^n_{{(Y,y)}}(X,x)\) is free for all \(x\in X\).
 \end{proof}

\begin{acknowledgement} We thank Rafael Holanda for his valuable comments on the manuscript, and Roger and Sylvia Wiegand for a suggestion given in a conversation about how to write some lines of the paper. Also, we are grateful to the referee for his/her careful reading of the manuscript and his/her appreciable suggestions.  All this contributed to the improvement of this paper. The first author was supported by grant 2022/03372-5, São Paulo Research Foundation (FAPESP). The second author was supported by grant 2019/21181-0, São Paulo Research Foundation (FAPESP). 
\end{acknowledgement}

\bibliographystyle{amsplain}
\bibliography{references.bib}
\end{document}